\documentclass[graybox]{svmult} 

\usepackage{helvet}         
\usepackage{courier}        
\usepackage{type1cm}        
\usepackage{makeidx}         
\usepackage{graphicx}        
\usepackage{multicol}        
\usepackage[bottom]{footmisc}

\usepackage{amsmath,amssymb,amsxtra,amscd,enumerate}
\usepackage[hypertex]{hyperref}   
\usepackage{longtable}
\usepackage{epic,eepic}
\usepackage{float}

\newtheorem{thm}{Theorem}[section]

\newtheorem{principle}[thm]{Guiding principle}
\newtheorem{exa}[thm]{Example}
\newtheorem{program}[thm]{Program}
\newtheorem{quest}[thm]{Question}
\newtheorem{probl}[thm]{Problem}
\makeatletter

\@addtoreset{equation}{section}
\makeatother

\begin{document}

\title*{Varna Lecture on $L^2$-Analysis of Minimal Representations}

\author{Toshiyuki Kobayashi}
\authorrunning{Kobayashi}
\institute{Toshiyuki Kobayashi \at Kavli IPMU and Graduate School of
Mathematical Sciences,
        The University of Tokyo, 3-8-1 Komaba, Meguro,
        153-8914 Tokyo, Japan}

\maketitle
\numberwithin{equation}{section}

\abstract{~Minimal representations of a real reductive group $G$ are
the \lq{smallest}\rq\ irreducible unitary representations of $G$. 
The author suggests a program
 of global analysis built
 on minimal representations from the philosophy:\\
   {\bf{small}} representation of a group = {\bf{large}}
symmetries in a representation space.  \\
This viewpoint serves
 as a driving force
 to interact algebraic representation theory
 with geometric analysis
 of minimal representations, 
yielding a rapid progress
 on the program.  
We give a brief guidance to recent works with emphasis
 on the Schr{\"o}dinger model.} 

\section{What are Minimal Representations?}
\label{sec:1}

Minimal representations
 of reductive groups $G$ are the \lq{smallest}\rq\
 infinite dimensional
 irreducible unitary representations.

The {\it{Weil}}
 (metaplectic,
 oscillator,
the Segal--Shale--Weil,
 harmonic)
 {\it{representation}},
known by a prominent role in number theory,
consists of two minimal representations of the metaplectic group
 $Mp(n,{\mathbb{R}})$.  
The minimal representation
 of a conformal group $SO(4,2)$
 arises on the
Hilbert space of bound states of the Hydrogen atom.

Minimal representations
 are distinguished 
 among other (continuously many) irreducible unitary representations
 of $G$
 by the following properties
 that I state loosely.
\begin{itemize}
\item[$\bullet$]
\lq{Smallest}\rq\ infinite dimensional representations of $G$.
\item[$\bullet$]
One of the \lq{building blocks}\rq\ of unitary representations
 of Lie groups.
\item[$\bullet$]
\lq{Closest}\rq\ to the trivial one dimensional representation
 of $G$.
\item[$\bullet$]
\lq{Quantization}\rq\
 of minimal nilpotent coadjoint orbits of $G$.
\item[$\bullet$]
Matrix coefficients have a \lq{slow decay}\rq\ at infinity.
\end{itemize}

In algebraic representation theory,
 there is a distinguished ideal ${\mathcal{J}}$
 introduced by Joseph \cite{xjoseph}
in the enveloping algebra 
 of a complex simple Lie algebra 
 other than type $A$
 (see also \cite{xgansavin}).  
An irreducible representation 
 of a real reductive Lie group $G$ 
 is called {\textit{minimal}}
 if its infinitesimal representation is annihilated
 by ${\mathcal{J}}$.  
Thus the terminology `minimal representations' is defined 
 {\it{inside}} representation theory.
We remark
 that not all reductive groups admit
 minimal representations.
Further,
 minimal representations
 are not always
 highest weight modules.
Beyond the case of highest weight modules,
 there has been an active study
 on minimal representations of reductive groups,
in particular,
 by algebraic approaches,
 see e.g., \cite{xgansavin,xjoseph, xKa,xkohcrcras,xKo,xMo,xTo,xV}.

In contrast,
my program focuses
 on global analysis
 inspired by minimal representations.  
For this,
 we switch the viewpoint,
 led by 

\begin{principle}
[\cite{xkgolden}]
\label{principle:1}
\begin{align*}
&\text{\textbf{small} representations of a group}
\nonumber
\\
={}
&\text{\textbf{large} symmetries in a representation space}.
\end{align*}
\end{principle}

An extremal case 
 of \lq{large symmetries}\rq\
 might be stated as
\begin{equation}
\label{eqn:dimmindim}
  \text{dimension of }\, \Xi
  <
  \text{dimension of 
 any non-trivial $G$-space}
\end{equation}
 when the representation of $G$
 is realized 
 on the space of functions on the geometry $\Xi$.  
An obvious implication of \eqref{eqn:dimmindim}
 is that $G$ cannot act on $\Xi$.  

The latter point of view,
served as a driving force,
 has brought us to a new line
 of investigation of geometric analysis modeled
 on minimal representations.  
In this program
 we are trying to dig out 
 new interactions
 with other areas of mathematics
 even {\textit{outside}} representation theory:

\begin{itemize}
\item[$\bullet$]
conformal geometry
 for general pseudo-Riemannian manifolds
 \cite{xkcheck, xkors1},
\item[$\bullet$]
Dolbeault cohomologies on open complex manifolds
 \cite{xkgqlim, xkohcrcras}.
\item[$\bullet$]
conservative quantities for PDEs
 \cite{xkcheck,xkors3},
\item[$\bullet$]
breaking symmetries and discrete branching laws
 \cite{xkors2, KOP, xkomin, xkowamo, xMo},
\item[$\bullet$]
 Schr\"odinger model and the unitary inversion operator
 \cite{xkhmL2, xkmano2, xkmanoAMS},
\item[$\bullet$]
deformation of the Fourier transform
 \cite{xbko},
\item[$\bullet$]
geometric quantization of nilpotent orbits
 \cite{xkhmL2,xkmanoAMS},
\item[$\bullet$]
holomorphic semigroup with a generalized Mehler kernel
 \cite{xbko, xkmano1, xkmano2},
\item[$\bullet$]
new orthogonal polynomials
 for fourth order differential operators
 \cite{xhkmm1, xhkmm2,xkmollers},
\item[$\bullet$]
a generalization of the Fock model and Bargmann transforms
 \cite{xhkmo}.
\end{itemize}

The aim of this article is
 to provide a brief guidance
 to the rapid progress
 on our program,
 \cite{xbko, xhkmm1,xhkmm2, xkhmL2,xhkmo,xkgqlim, xkgolden,xkmanoAMS, xkmollers, KOP, xkowamo}.
We should mention
 that in order to avoid an overlap
 with a recent publication \cite{xkgolden},
 we do not include here some other constructions
 such as a conformal model
 of minimal representations
 (e.g. the construction of the intrinsic
 conservative quantities
 for the conformally invariant differential equations).
Instead,
 we highlight an $L^2$-model
 ({\textit{Schr{\"o}dinger model}})
 of the minimal representations
 and its variant.  
We apologize 
 for not being able to mention
 some other important works
 on minimal representations,
 e.g., 
 see \cite{xgansavin}
 and references therein.  
For a comparison 
 of the $L^2$-model 
 with the conformal model,
 we refer to 
 \cite[Introduction]{xkmanoAMS}.

\section{More Symmetric than Symmetric Spaces}
\label{sec:moresym}

The traditional geometric construction
 of representations
 of Lie groups $G$ is given
 in the following two steps:
\newline
Step 1.\enspace
The group $G$ acts on a geometry $X$.
\newline
Step 2.\enspace
By the translation,
 $G$ acts linearly on the space $\Gamma(X)$
 of functions
 (sections of equivariant bundles,
 or cohomologies, 
 $\cdots$).

Na{\"i}vely,
 the Gelfand--Kirillov dimension
 of the representation on $\Gamma(X)$
 is supposed to be 
 the dimension of $X$.
Thus we may expect
 that the representation on the function space $\Gamma(X)$
 is \lq{small}\rq\
 if the geometry $X$ itself is small.

First of all,
we ask when the geometry $X$ is \lq{small}\rq.

For this we may begin with the case when $G$ acts transitively on $X$,
 or equivalently,
 $X$ is a homogeneous space $G/H$.
Further,
 if we compare two homogeneous spaces
 $X_1=G/H_1$ and $X_2=G/H_2$
 with $H_1 \subset H_2$,
 we may think that $X_2$ is smaller than $X_1$.
Hence \lq{smaller}\rq\ representations on $\Gamma(X)$
 should be attained
 if $X=G/H$
 where $H$ is a maximal subgroup of $G$.

Here are two typical settings
 for real reductive Lie groups $G$:
\begin{itemize}
\item[$\bullet$]
$(G,H)$ is a symmetric pair.

In this case,
 the Lie algebra ${\mathfrak {h}}$
 of $H$ is maximal reductive in ${\mathfrak {g}}$.
Analysis on reductive symmetric spaces $G/H$
 has been largely developed in
 particular,
 since 1950s
 by the Gelfand school,
 Harish-Chandra,
 Shintani,
 Helgason,
 Takahashi, 
 Molchanov, 
Faraut, 
Flensted-Jensen,
 Matsuki--Oshima--Sekiguchi,
 Delorme,
 van den Ban,
Schlichtkrull,
 among others.
\item[$\bullet$]
$H$ is a Levi subgroup of $G$.

In this case,
 there exists a $G$-invariant polarization on $G/H$,
 and its geometric quantization
 obtained by the combination of the Mackey induction
 (real polarization)
 and the Dolbeault cohomologies
 (complex polarization)
 produces a \lq{generic part}\rq\
 of irreducible unitary representations of $G$.
The resulting representations
 are the \lq{smallest}\rq\
 if $H$ is a maximal Levi subgroup.
\end{itemize}
These two typical examples
 are related:
Tempered representations
 for reductive symmetric spaces
 (i.e.~irreducible unitary representations
 that contribute to $L^2(G/H)$)
 are given by the combination
 of the ordinary and cohomological 
 parabolic inductions.  
A missing picture 
 in the above two settings
 is so called 
 \lq{unipotent representations}\rq\
 including minimal representations.  

On the other hand,
 it is rare but still happens
 that the representation of $G$ on the function space
 $\Gamma(X)$ extends
 to a representation of a group $\widetilde G$
 which contains $G$,
 even when the $G$-action on the geometry $X$
 does not extend to $\widetilde G$
(in particular, 
 Step 1 does not work
 for the whole group $\widetilde G$).  
We discuss this phenomenon
 in the Schr{\"o}dinger model
 of minimal representations
 when $G$ is a maximal parabolic subgroup
 (the notation $(G, \widetilde G)$
 here will be replaced
 by $(P,G)$ in Section \ref{sec:schrodinger}).
Such a phenomenon also occurs
 when $G$ is reductive.  
Thus the analysis of minimal representations
 may be thought of
 as \lq{analysis
 with more symmetries}\rq\
 than the traditional analysis
 on homogeneous spaces.
Here is a typical example:
\begin{exa}
[{\cite[Theorem 5.3]{xkcheck}}]
The minimal representation of the indefinite orthogonal group
 $\widetilde G=O(p,q)$
 ($p+q$:even)
 is realized in function spaces
 on symmetric spaces
 of the subgroups $G=O(p-1,q)$ or $O(p,q-1)$
 on which the whole group $\widetilde G$
 cannot act geometrically.
\end{exa}

\begin{exa}
[\cite{KOP}]
The restriction of the most degenerate principal series representations
of $\widetilde G=GL(n,{\mathbb{R}})$
 to the subgroup
$G=O(p,q)$
 ($p+q=n$)
 reduces to the analysis
 of the symmetric space of $G$
 on which the whole group $\widetilde G$
 cannot act transitively.
\end{exa}

Further examples and explicit branching rules
 can be found in \cite{xkcheck, xkors2, KOP}
 where the restriction of minimal representations
 to subgroups
 ({\it{broken symmetries}})
 reduce to analysis
 on certain semisimple symmetric spaces.

\section{Schr\"odinger Model of Minimal Representations}
\label{sec:schrodinger}

Any coadjoint orbit
 of a Lie group is naturally 
 a symplectic manifold
 endowed with the Kirillov--Kostant--Souriau
 symplectic form.
For a reductive Lie group $G$, 
\lq{geometric quantization}\rq\
of semisimple coadjoint orbits has been considerably well-understood
 --- 
this corresponds to the ordinary
 or cohomological parabolic induction
 in representation theory,
 whereas \lq{geometric quantization}\rq\
 of nilpotent coadjoint orbits
 is more mysterious 
 (see \cite{xBK, xhkmo, xkgqlim}).

In this section
 we explain a recent work 
 \cite{xkhmL2}
 with Hilgert and M{\"o}llers
 on the $L^2$-construction
of minimal representations
 built on a Lagrangian subvariety
 of a real minimal nilpotent orbit, 
 which continues 
 a part of the earlier works
  \cite{xkors3} with {\O}rsted, 
 and \cite{xkmanoAMS} with Mano.

Suppose that $V$ is a simple Jordan algebra
 over ${\mathbb{R}}$.
We assume
 that its maximal Euclidean Jordan subalgebra
is also simple.
Let $G$ and $L$ be the identity components
 of the conformal group
 and the structure group of the Jordan algebra $V$,
 respectively.
Then the Lie algebra ${\mathfrak {g}}$
 is a real simple Lie algebra
 and has a Gelfand--Naimark decomposition
$ {\mathfrak {g}}= \overline{{\mathfrak {n}}}
                    + {\mathfrak {l}}
                    + {\mathfrak {n}},
$
where ${\mathfrak {n}} \simeq V$
 is regarded as an Abelian Lie algebra,
 ${\mathfrak {l}} \simeq {\mathfrak {str}}(V)$
 the structure algebra,
 and $\overline{{\mathfrak {n}}}$
 acts on $V$
 by quadratic vector fields.

Let ${\mathbb{O}}_{\operatorname{min}}^{G_{\mathbb{R}}}$
 be a (real) minimal nilpotent coadjoint orbit.
By identifying ${\mathfrak {g}}$
 with the dual ${\mathfrak {g}}^{\ast}$,
 we consider the intersection 
 $V \cap {\mathbb{O}}_{\operatorname{min}}^{G_{\mathbb{R}}}$, 
 which may be disconnected 
 (this happens in the case \eqref{eqn:J2} below).  
Let $\Xi$ be any connected component
 of 
$
  V \cap {\mathbb{O}}_{\operatorname{min}}^{G_{\mathbb{R}}}.
$
We note that the group $L$ acts on $\Xi$
 but $G$ does not.  
There is a natural $L$-invariant Radon measure
 on $\Xi$, 
 and we write $L^2(\Xi)$
 for the Hilbert space 
 consisting of square integrable functions 
 on $\Xi$.  
Then we can define 
 a unitary representation
 on $L^2(\Xi)$
 (Schr\"odinger model)
 built on a Lagrangian submanifold $\Xi$
 in this generality \cite{xkhmL2},
 see also \cite{xDS, xkors3}.  
\begin{thm}
[Schr\"odinger model]
\label{thm:khmL2}
Suppose $V \not \simeq {\mathbb{R}}^{p,q}$
 with $p+q$ odd.
\par\noindent
{\rm{1)}}\enspace
$\Xi$ is a Lagrangian submanifold
 of ${\mathbb{O}}_{\operatorname{min}}^{G_{\mathbb{R}}}$.
\par\noindent
{\rm{2)}}\enspace
There is a finite covering group $G\,\widetilde{}$
 of $G$
 such that $G\,\widetilde{}$ acts on $L^2(\Xi)$
 as an irreducible unitary representation.  
\par\noindent
{\rm{3)}}\enspace
The Gelfand--Kirillov dimension of $\pi$
attains its minimum
among all infinite dimensional representations
 of $G\,\widetilde{}$, 
 i.e. $\operatorname{DIM}(\pi)=\frac 1 2 \dim {\mathbb{O}}_{\operatorname{min}}^{G_{\mathbb{R}}}$.  
\par\noindent
{\rm{4)}}\enspace
The annihilator of the differential representation
 $d \pi$ is the Joseph ideal
 in the enveloping algebra $U({\mathfrak {g}})$
 if $V$ is split and ${\mathfrak {g}}$ is not of type $A$.
\end{thm}

The simple Lie algebras
 arisen in Theorem \ref{thm:khmL2}
 are listed as follows:
\begin{align}
&{\mathfrak {sl}}(2k,{\mathbb{R}}),
{\mathfrak {so}}(2k,2k),
{\mathfrak {so}}(p+1,q+1),
{\mathfrak {e}}_{7(7)},
\label{eqn:J1}
\\
&
{\mathfrak {sp}}(k,{\mathbb{R}}),
{\mathfrak {su}}(k,k),
{\mathfrak {so}}^{\ast}(4k),
{\mathfrak {so}}(2,k),
{\mathfrak {e}}_{7(-25)},
\label{eqn:J2}
\\
&{\mathfrak {sp}}(k,{\mathbb{C}}),
{\mathfrak {sl}}(2k,{\mathbb{C}}),
{\mathfrak {so}}(4k,{\mathbb{C}}),
{\mathfrak {so}}(k+2,{\mathbb{C}}),
{\mathfrak {e}}_{7}({\mathbb{C}}),
\label{eqn:J3}
\\
&{\mathfrak {sp}}(k,k),
{\mathfrak {su}}^{\ast}(4k),
{\mathfrak {so}}(k,1).
\label{eqn:J4}
\end{align}

\begin{remark}
In the case where $V$ is an Euclidean Jordan algebra, 
 $G$ is the automorphism group
 of a Hermitian symmetric space
 of tube type
 (see \eqref{eqn:J2})
 and there are two real minimal nilpotent orbits.
The resulting representations $\pi$
are highest (or lowest) weight modules.
\end{remark}
\begin{remark}
If the complex minimal nilpotent orbit
${\mathbb{O}}_{\operatorname{min}}^{G_{\mathbb{C}}}$
 intersects with ${\mathfrak {g}}$,
 then ${\mathbb{O}}_{\operatorname{min}}^{G_{\mathbb{R}}}$
 is equal to
$
  {\mathbb{O}}_{\operatorname{min}}^{G_{\mathbb{C}}} \cap {\mathfrak {g}}
$
 or its connected component.
We notice that
$
  {\mathbb{O}}_{\operatorname{min}}^{G_{\mathbb{C}}} \cap {\mathfrak {g}}
$
 may be an empty set
 depending on the real form ${\mathfrak {g}}$.  
In the setting of Theorem \ref{thm:khmL2},
this occurs
 for \eqref{eqn:J4}.
In this case,
 the representation $\pi$ 
 in Theorem \ref{thm:khmL2}
 is not a minimal representation
 as the annihilator of $d \pi$ is not
 the Joseph ideal,
 but $\pi$ is still one of the \lq{smallest}\rq\
 infinite dimensional representations
 in the sense
 that the Gelfand--Kirillov dimension
 attains its minimum.
\end{remark}

\begin{remark}
There is no minimal representation
 for any group
 with Lie algebra ${\mathfrak {o}}(p+1,q+1)$
 with $p+q$ odd,
 $p$, $q$ $\ge 3$
 (see \cite[Theorem 2.13]{xV}).
\end{remark}

\begin{exa}
\label{exam:symmsp}
Let $V= \operatorname{Sym}(m,{\mathbb{R}})$.
Then
$G=Sp(m,{\mathbb{R}})$
 and
\begin{equation}
\label{eqn:coneSp}
V \cap {\mathbb{O}}_{\operatorname{min}}^{G_{\mathbb{R}}}
=\{X \in M(m,{\mathbb{R}}):
 X = {}^{t\!}X, \operatorname{rank}X=1\}.
\end{equation}
Let 
$
   \Xi
:=\{X \in V \cap {\mathbb{O}}_{\operatorname{min}}^{G_{\mathbb{R}}}:
    \operatorname{Trace}X>0\}.
$
Then the double covering map ({\it{folding map}})
\[
  {\mathbb{R}}^m \setminus\{0\} \to \Xi,
  \quad
  v \mapsto v {}^{t\!}v
\]
 induces an isomorphism between $L^2(\Xi)$ and the Hilbert space
 $L^2({\mathbb{R}}^m)_{\operatorname{even}}$
 of even square integrable functions
 on ${\mathbb{R}}^m$.
Thus our representation $\pi$
 on $L^2(\Xi)$ is nothing
 but the Schr{\"o}dinger model of the even part
 of the Segal--Shale--Weil representation
 of the metaplectic group $Mp(m,{\mathbb{R}})$
 \cite{xfolland, Howe}.
\end{exa}
\begin{exa}
\label{exam:4.2}
Let $V= {\mathbb{R}}^{p,q}$
 with $p+q$ even.
Then
 ${\mathfrak {g}}={\mathfrak {o}}(p+1,q+1)$,
 and
\begin{equation}
\label{eqn:coneO}
     V \cap {\mathbb{O}}_{\operatorname{min}}^{G_{\mathbb{R}}}
=\{\xi \in {\mathbb{R}}^{p+q}:
 \xi_1^2+ \cdots +\xi_p^2 -\xi_{p+1}^2 - \cdots - \xi_{p+q}^2 =0\}
\setminus \{0\}.
\end{equation}
If $p=1$, 
 $V \cap {\mathbb{O}}_{\operatorname{min}}^{G_{\mathbb{R}}}$
 consists of two connected components
 according to the signature
of $\xi_1$,
 i.e.~the past and future cones.  
They yield highest/lowest weight modules.  
For $p,q \ge 2$, 
 $V \cap {\mathbb{O}}_{\operatorname{min}}^{G_{\mathbb{R}}}$
 is connected,
 and our representation $\pi$ on $L^2(\Xi)$ is 
 the Schr{\"o}dinger model
 of the minimal representation
 of $O(p+1,q+1)$
 constructed in \cite{xkors3}, 
 which is a neither highest nor lowest weight module.
\end{exa}

As we discussed in Section \ref{sec:moresym}
in contrast to traditional analysis
 on homogeneous spaces,
 the group $G$ in our setting
 is too large to act geometrically on $\Xi$.
This very feature
 in the Schr{\"o}dinger model
 is illustrated by the fact
 that the Lie algebra $\overline{\mathfrak {n}}$
 acts as differential operators on $\Xi$
 of second order.
They are {\it{fundamental differential operators}}
 \cite{xkmanoAMS}
 in the setting of Example \ref{exam:4.2}
 (see also Bargmann--Todorov
 \cite{xbt}).
In \cite{xkhmL2},
 these differential operators are said to be 
{\it{Bessel operators}},
 and serve as a basic tool
 to study the Schr{\"o}dinger model $\pi$
 in the setting of Theorem \ref{thm:khmL2}.

\section{Special Functions to 4th order Differential Operators}

Guiding Principle \ref{principle:1}
 suggests that there should exist
 plentiful functional equations
 in the representation spaces
 for minimal representations.
Classically,
 it is well-known
 that Hermite polynomials form
 an orthogonal basis
 for the radial part of the Schr{\"o}dinger model
 of the Weil representation \cite{xfolland},
 whereas Laguerre polynomials arise in the minimal representation
 of the conformal group $SO(n,2)$
 (\cite{xtodorov}).

These classical minimal representations
 are highest weight modules.
However,
 for more general reductive groups,
 minimal representations
 do not always have highest weight vectors, 
 and the corresponding \lq{special functions}\rq\
 do not necessarily satisfy
 second order differential equations.
We found in \cite{xkmanoAMS}
that Meijer's $G$-functions
 $G_{0\,4}^{2\,0}(x|b_1,b_2,b_3,b_4)$
 play an analogous role in the minimal representation
 of $O(p,q)$.  
Here Meijer's $G$-functions $G_{0\,4}^{2\,0}(x|b_1, b_2, b_3, b_4)$
satisfy a fourth order ordinary differential equation
\[
     \prod_{j=1}^{4} (x\frac{d}{dx}-b_j)u(x)=xu(x).  
\]

More generally,
 the following fourth order differential operators 
\[
  {\mathcal{D}}_{\mu,\nu}
  :=
  \frac{1}{x^2}
  ((\theta + \nu)(\theta + \mu + \nu)-x^2)
  (\theta(\theta + \mu)-x^2)
  -\frac{(\mu-\nu)(\mu+ \nu +2)}{2}
\]
appear
 naturally
 in the Schr{\"o}dinger model
 of minimal representations
 in the setting of Theorem \ref{thm:khmL2}.
Here $\theta = x \frac {d}{dx}$.

The subject of \cite{xhkmm1, xhkmm2, xkmollers}
 is the study of eigenfunctions of ${\mathcal{D}}_{\mu,\nu}$
 including
\begin{itemize}
\item[$\bullet$]
generating functions
 for eigenfunctions of ${\mathcal{D}}_{\mu,\nu}$,
\item[$\bullet$]
asymptotic behavior near the singularities,

\item[$\bullet$]
$L^2$-eigenfunctions
 and concrete formulas of $L^2$-norms,

\item[$\bullet$]
integral representations
 of eigenfunctions,

\item[$\bullet$]
recurrence relations among eigenfunctions,
\item[$\bullet$]
(local) monodromy.
\end{itemize}
The $L^2$-eigenfunctions of ${\mathcal{D}}_{\mu,\nu}$
 arise as $K$-finite vectors
 in the Schr{\"o}dinger model of the minimal representations
 constructed in Theorem \ref{thm:khmL2}
 in a uniform fashion.
These \lq{special functions}\rq\
 with certain integral parameters
 yield orthogonal polynomials
 (the {\it{Mano polynomials}} $M_{j}^{\mu,l}(x)$)
satisfying again fourth order differential equations
 \cite{xhkmm2},
 which include Hermite polynomials
 and Laguerre polynomials
 as special cases.
We note
 that the fourth order differential equation
$
  {\mathcal{D}}_{\lambda,\mu}f=\nu f
$
 reduces to a differential equation 
 of second order
 when $G/K$ is a tube domain
 (see \eqref{eqn:J2}).  
See also Kowata--Moriwaki \cite{xkowamo}
 for further analysis 
 of the fundamental differential operators
 on $\Xi$.

\section{Broken Symmetries and Branching Laws}
As indicated in Guiding Principle \ref{principle:1}, 
the `large symmetries' in representation spaces
 of minimal representations
 produce also fruitful examples
 of branching laws
 which we can expect a simple
 and detailed study.

Suppose $\pi$
 is a unitary representation
 of a real reductive Lie group $G$.
We consider $\pi$ as a representation
 of a subgroup $G'$ of $G$, 
 referring it as the restriction $\pi|_{G'}$.  
In general,
 the restriction $\pi|_{G'}$ decomposes into a direct integral of
 irreducible representations
 of $G'$
 ({\textit{branching law}}).  
It often happens that the branching law
 contains continuous spectrum
 if $G'$ is non-compact.
Even worse,
 each irreducible representation
 of $G'$ may occur in the branching law with infinite multiplicities.
See \cite{xkaspm} for such wild examples
 even when $(G,G')$ is a symmetric pair.
In \cite{xk:1, xk:2},
 we raised the following:
\begin{program}
\label{prog:deco}
1) \enspace
Determine the triple $(G, G',\pi)$
 for which the restriction $\pi|_{G'}$ decompose discretely
 with finite multiplicities.
\par\noindent
2)\enspace
 Find branching laws
 for (1).
\end{program}

Program \ref{prog:deco}
 intends to single out a nice framework of branching problems
 for which we can expect
 a detailed
 and explicit study of the restriction.
Concerning Program \ref{prog:deco} (1)
 for Zuckerman's derived functor modules $\pi$,
 a necessary and sufficient condition
 for {\it{discrete decomposition
 with finite multiplicities}}
was proved 
in \cite{xk:2, xk:4},
 and a complete classification
 was given with Oshima \cite{xkoadv}
 when $(G,G')$ is a reductive symmetric pair.  

As such,
 the local theta correspondence
 with respect to compact dual pairs
 is a classic example
 for minimal representations $\pi$:

\begin{exa}
Suppose that $\pi$ is the Weil representation,
 and that $G'=G_1' \cdot G_2'$ is a dual pair
 in $G=Mp(n,{\mathbb{R}})$
 with $G_2'$ compact.
Then the restriction $\pi|_{G'}$
 decomposes discretely
 and multiplicity-freely.
The resulting branching laws yield
 a large part of unitarizable highest weight modules
 of $G_1'$ (Enright--Howe--Wallach \cite{xEHW}).
\end{exa}
In order to discuss Program \ref{prog:deco}
 for minimal representations,
 we recall from \cite{xk:2, xk:3, xk:4} the general theory.  
Let $K$ be a maximal compact subgroup of $G$,
 $T$ a maximal torus of $K$,
 and ${\mathfrak {t}}$, ${\mathfrak {k}}$ the Lie algebras
 of $T$, $K$,
 respectively.
We choose the set $\Delta^+({\mathfrak {k}}, {\mathfrak {t}})$
 of positive roots,
 and denote by ${\mathfrak {t}}_+$
 the dominant Weyl chamber in $\sqrt{-1}{\mathfrak {t}}^*$.
We also fix a $K$-invariant inner product on ${\mathfrak {k}}$,
 and regard $\sqrt{-1} {\mathfrak {t}}^*$ as a subspace
 of $\sqrt{-1}{\mathfrak {k}}^*$.

First,
 suppose that $K'$ is a closed subgroup of $K$.
The group $K$ acts on the cotangent bundle $T^{\ast}(K/K')$
 of the homogeneous space $K/K'$
 in a Hamiltonian fashion.
 We write
$$
  \mu: T^{\ast}(K/K') \to \sqrt{-1}{\mathfrak {k}}^*
$$
for the momentum map,
 and define the following closed cone by
$$
    C_K(K') := \operatorname{Image}\mu \cap {\mathfrak {t}}_+.
$$

Second,
let $\operatorname{Supp}_K(\pi)$ be the set
 of highest weights of finite dimensional irreducible representations
 of $K$ occurring in a $K$-module $\pi$.
The asymptotic $K$-support ${\operatorname{AS}}_K(\pi)$
 is defined to be the asymptotic cone
 of $\operatorname{Supp}_K(\pi)$.
It is a closed cone in ${\mathfrak {t}}_+$.
There are only finitely many possibilities
 of ${\operatorname{AS}}_K(\pi)$
 for the restriction $\pi|_K$
 of irreducible representations $\pi$ of $G$.

The asymptotic cone ${\operatorname{AS}}_K(\pi)$
 tends to be a \lq{small}\rq\ subset
in ${\mathfrak {t}}_+$
 if $\pi$ is a \lq{small}\rq\ representation.
For example,
\begin{alignat}{2}
{\operatorname{AS}}_K(\pi)
  =&
  \{0\}
  &&\text{ if } \dim \pi < \infty,
\notag
\\
{\operatorname{AS}}_K(\pi)
  =&
  {\mathbb{R}}_+ \beta
  \qquad
  &&\text{ if $\pi$ is a minimal representation,}
\label{eqn:AS}
\end{alignat}
where $\beta$ is the highest root
 of the $K$-module
 ${\mathfrak {p}}_{\mathbb{C}}
  :={\mathfrak {g}}_{\mathbb{C}}/{\mathfrak {k}}_{\mathbb{C}}
$.
The formula \eqref{eqn:AS}
 holds in a slightly more general setting
 where the associated variety
 of $\pi$
 is the closure of a single minimal nilpotent $K_{\mathbb{C}}$-orbit
 on ${\mathfrak{p}}_{\mathbb{C}}$
 \cite{xkomin}.
Concerning Program \ref{prog:deco},
 we established an easy-to-check
 criterion in \cite{xk:3}:

\begin{thm}
\label{thm:deco}
Suppose $G'$ is a reductive subgroup of $G$
 such that $K' := G' \cap K$
 is a maximal compact subgroup of $G'$.
If
\begin{equation}
\label{eqn:criterion}
C_K(K') \cap {\operatorname{
AS}}_K(\pi) = \{0\},
\end{equation}
then the restriction $\pi|_{G'}$ decomposes discretely into a direct sum
of irreducible unitary representations of $G'$
 with finite multiplicities.
\end{thm}

As was observed in \cite{xkprog05}, 
 we can expect from  the formula \eqref{eqn:AS}
 and from the criterion \eqref{eqn:criterion}
 that there is plenty of subgroups $G'$ for which the restriction of the minimal representation of $G$ decomposes discretely
 and with finite multiplicities.
Reductive symmetric pairs
 $(G,G')$
 for which the restriction $\pi|_{G'}$
 is (infinitesimally) discretely decomposable
 for a minimal representation $\pi$ of $G$
 has been recently classified in \cite{xkomin}.

\section{Generalized Fourier Transform as a Unitary Inversion}
\label{sec:genfourier}

In the $L^2$-model of the minimal representation $\pi$
 of $G$ on
$L^2(\Xi)$,
the action of the maximal parabolic subgroup $P$
 with Lie algebra ${\mathfrak {l}}+{\mathfrak {n}}$ is simple,
namely,
it is given just by translations and multiplications.
Let $w$ be  the conformal inversion
 of the Jordan algebra.
In light of the Bruhat decomposition
\[G = P \amalg PwP,
\]
 it is enough to find $\pi(w)$
in order to give a global formula
 of the $G$-action on $L^2(\Xi)$.
We highlight this specific unitary operator,
 and set
\begin{equation}\label{eqn:FXi}
\mathcal{F}_\Xi := c\pi(w),
\end{equation}
where $c$ is a complex number of modulus one
 (the {\it{phase factor}}).
We call ${\mathcal{F}}_\Xi$
 the {\it{unitary inversion operator}}.
We studied
 in a series of papers \cite{xkmano1, xkmano2, xkmanoAMS}
 with Mano the following:
\begin{probl}
\label{prog:6.1}
Find an explicit formula of the integral kernel of ${\mathcal{F}}_\Xi$.
\end{probl}

The kernel of the Euclidean Fourier transform
 is given by $e^{-i\langle x, \xi \rangle}$, 
 which is locally integrable.  
It is plausible
 that this analytic feature happens 
 if and only if 
 the corresponding minimal representation
 is of highest weight.  
Thus we raise the following:

\begin{quest}
Let $(\pi, L^2(\Xi))$ be the $L^2$-model
 of a minimal representation $\pi$ 
 of a simple Lie group $G\,\widetilde{}$
 constructed on a Lagrangian submanifold $\Xi$
 of $\mathbb{O}_{\operatorname{min}}^G$
 as in Theorem \ref{thm:khmL2} {\rm{\cite{xkhmL2}}}.  
Are the following two conditions equivalent?
\begin{enumerate}
\item[{\rm{(i)}}]
The kernel of the unitary inversion operator
 ${\mathcal{F}}_{\Xi}$ is locally integrable.  
\item[{\rm{(ii)}}]
$\pi$ is a highest/lowest weight module.  
\end{enumerate}
\end{quest}
Here we have excluded the case
 where the simple Lie algebra ${\mathfrak{g}}$
 is of type $A_n$
 (the Joseph ideal is not defined 
 for ${\mathfrak {g}}_{\mathbb{C}}={\mathfrak {sl}}_n({\mathbb{C}})$).  
In the case $G=O(p+1, q+1)$
 with $p+q$ even $>2$, 
 it was proved in \cite{xkmanoAMS}
 that (i) holds
 if and only if either $\operatorname{min}(p,q)=1$
 (equivalently, (ii) holds)
 or $(p,q)=(3,3)$ 
 (equivalently,
 ${\mathfrak {g}}={\mathfrak {o}}(3,3)
  \simeq {\mathfrak {sl}}(4,{\mathbb{R}})$
 is of type $A_3$).  
The implication (ii) $\Rightarrow$ (i) was 
 proved in \cite{xhkmo} for tube type,  
see \eqref{eqn:Ftube}.  
The implication (i) $\Rightarrow$ (ii) is an open problem
 except for the above mentioned case
 ${\mathfrak{g}} ={\mathfrak {o}}(p+1,q+1)$.  

When $G=O(p+1,q+1)$
(see Example \ref{exam:4.2}),
$\mathcal{F}_\Xi$ intertwines the multiplication of
coordinate functions $\xi_j$ $(1 \le j \le p+q)$
 with the operators $R_j$ $(1 \le j \le p+q)$ which are mutually commuting
differential operators of second order on $\Xi$
(see Bargmann--Todorov \cite{xbt}, \cite[Chapter 1]{xkmanoAMS}).

This algebraic feature is similar to the classical fact
 that the Euclidean Fourier transform
$\mathcal{F}_{\mathbb{R}^m}$ intertwines the multiplication operators
$\xi_j$ and the differential operators $\sqrt{-1}\partial_j$
$(1 \le j \le m)$
 (see Example \ref{exam:symmsp}).
In the setting of Theorem \ref{thm:khmL2},
 ${\mathcal{F}}_{\Xi}$ intertwines
 the multiplication
 of coordinate functions
 with Bessel operators.
Actually,
 this algebraic feature determines uniquely 
 ${\mathcal{F}}_{\Xi}$ 
 up to a scalar
 \cite{xkhmL2, xkmanoAMS}.

Concerning Problem \ref{prog:6.1}, 
the first case is well-known
 (see \cite{xfolland} for example):
\par\noindent
{1)}\enspace
${\mathfrak {g}}={\mathfrak {sp}}(m,{\mathbb{R}})$.

${\mathcal{F}}_\Xi
= \text{the Euclidean Fourier transform on ${\mathbb{R}}^m$}$.

Here are some recent results
 on a closed formula
 of the integral kernel:
\par\noindent
{2)}\enspace
${\mathfrak {g}}={\mathfrak {o}}(p+1, q+1)$
  (with Mano \cite{xkmano2}).
\par\noindent
{3)}\enspace
The associated Riemannian symmetric space $G/K$
 is of tube type
 (see \eqref{eqn:Ftube}).

We note that
 minimal representations in the cases 1) and 3)
 are highest (or lowest) weight modules,
 whereas minimal representations
 in the case 2) do not have highest weights
 when $p,q \ge 2$ and $p+q$ is odd.

Problem \ref{prog:6.1} is open for other cases,
 in particular, for minimal representations
 without highest weights
 except for the case $G=O(p+1,q+1)$.

\section{$SL_2$-triple in the Schr{\"o}dinger Model}
\label{sec:sl2tri}
On ${\mathbb{R}}^m$, 
 we set $|x|:=(\sum_{j=1}^{m}x_j^2)^{\frac1 2}$, 
 $E:=\sum_{j=1}^{m} x_j \frac{\partial}{\partial x_j}$
 (Euler operator)
 and $\Delta = \sum_{j=1}^{m} \frac{\partial^2}{\partial x_j^2}$
(Laplacian).  
Then 
it is classically known
 (e.g., \cite{xfolland,Howe}) that the operators
\begin{equation}
\label{eqn:Spsl}
\widetilde h':= E+\frac{m}2,
\quad
\widetilde e':=\frac{\sqrt{-1}}2|x|^2,
\quad
\widetilde f':=\frac{\sqrt{-1}}2 \Delta
\end{equation}
form an ${\mathfrak {sl}}_2$-triple, 
namely,
the following commutation relation holds:
\[
  [\widetilde{h}', \widetilde{e}']=2\widetilde{e}',
  \quad 
  [\widetilde{h}', \widetilde{f}']=-2\widetilde{f}', 
  \quad
  [\widetilde{e}', \widetilde{f}']=\widetilde{h}'.  
\]

On the other hand,
 we showed in \cite{xkmano2}
 that the following operators
\begin{equation}
\label{eqn:Osl}
  \widetilde h:=2 E +m-1,
\quad
  \widetilde e:=2{\sqrt{-1}}|x|,
 \quad
  \widetilde f:=\frac{\sqrt{-1}}2 |x| \Delta
\end{equation}
also forms an ${\mathfrak {sl}}_2$-triple, 
i.e., 
$
[\widetilde h, \widetilde e]=2 \widetilde e,
\quad
[\widetilde h, \widetilde f]= -2 \widetilde f,
\quad
[\widetilde e, \widetilde f]=\widetilde h.  
$

Further the differential operator
\[
D:=\frac{1}{2\sqrt{-1}}(-\widetilde e+\widetilde f)
 =|x|(\frac{\Delta}{4}-1)
\]
extends to a self-adjoint operator
 and has only discrete spectra
 on $L^2({\mathbb{R}}^m, \frac{dx}{|x|})$
 which are given by
$
\{-(j+\frac {m-1}{2}): j=0,1,2, \cdots\}
$
 (see \cite{xkmano2}), 
whereas the {\it{Hermite operator}}
\[
 \mathcal D:=\frac{1}{2\sqrt{-1}}(-\widetilde e'+\widetilde f')
            =\frac 1 4(\Delta-|x|^2)
\]
extends to a self-adjoint operator
 and has only discrete spectra
 on $L^2({\mathbb{R}}^m, dx)$
 which are given by
$
\{-\frac 1 2(j+\frac {m}{2}): j=0,1,2, \cdots\}
$
(see \cite{xfolland, Howe}).  
Hence,
 one can define
 for $\operatorname{Re}t \ge 0$:
\begin{align*}
e^{tD} :=& \sum_{k=0}^\infty \frac{t^k}{k!} D^k
\quad
\text{on } L^2({\mathbb{R}}^m, \frac{dx}{|x|}),
\\
e^{t{\mathcal{D}}} :=& \sum_{k=0}^\infty \frac{t^k}{k!} {\mathcal{D}}^k
\quad
\text{on } L^2({\mathbb{R}}^m, dx).  
\end{align*}
They are holomorphic one-parameter semigroups
 consisting of Hilbert--Schmidt operators
 for $\operatorname{Re}t >0$,
 and are unitary operators for $\operatorname{Re}t =0$.

A closed formula for both $e^{t D}$
 and $e^{t {\mathcal{D}}}$
 is known.  
That is, 
the holomorphic semigroup $e^{t{\mathcal{D}}}$
 has the classical Mehler kernel
 given by the Gaussian kernel $e^{-|x|^2}$
 and reduces to the Euclidean Fourier transform
 when $t=\sqrt{-1}\pi$
 (\cite[\S 5]{Howe}
),
whereas the integral kernel of the holomorphic semigroup $e^{tD}$
 is given by the $I$-Bessel function
 and the special value
 at $t=\sqrt{-1}\pi$
 is by the $J$-Bessel function
 (see \cite[Theorem A and Corollary B]{xkmano2}
 for concrete formulas).

We can study these holomorphic semigroups
 by using the theory of discretely decomposable
 unitary representations
 (e.~g.~\cite{xk:1,xk:2,xk:3}).
Actually,
the aforementioned ${\mathfrak {sl}}_2$-triple
 arises as the differential action of
 the Schr{\"o}dinger model
 of the minimal representations
 of $Mp(m,{\mathbb{R}})$ on $L^2({\mathbb{R}}^m,dx)$
 and $SO_0(m+1,2)$ on $L^2({\mathbb{R}}^m,\frac{dx}{|x|})$,
respectively
 via
\begin{alignat*}{2}
  {\mathfrak {sl}} (2,{\mathbb{R}})
  &\simeq {\mathfrak {sp}} (1,{\mathbb{R}})
  && \subset {\mathfrak {sp}}(m,{\mathbb{R}}),
\\
  {\mathfrak {sl}} (2,{\mathbb{R}})
  & \simeq {\mathfrak {so}}(1,2)
  && \subset {\mathfrak {so}}(m+1,2),
\end{alignat*}
for 
which we write as 
 $d \iota : {\mathfrak {sl}}(2,{\mathbb{R}}) \hookrightarrow {\mathfrak {g}}$.

In both cases,
 the Lie algebra ${\mathfrak {g}}$ contains
 a subalgebra commuting with $\iota({\mathfrak {sl}} (2,{\mathbb{R}}))$, 
 which is isomorphic to ${\mathfrak {o}}(m)$.    
Then the minimal representations
 decompose as the representation
 of the direct product group $SL (2,{\mathbb{R}}) \times O(m)$
 (up to coverings and connected groups) as follows:

\begin{align*}
  L^2(\mathbb R^m, \frac{dx}{|x|}) \simeq
  \sideset{}{^\oplus}{\sum}_{j=0}^{\infty}
  \pi_{2j+m-1}^{SL(2, \mathbb R)} \boxtimes \mathcal H^{j}({\mathbb{R}}^m).
\\
  L^2(\mathbb R^m, dx) \simeq
  \sideset{}{^\oplus}{\sum}_{j=0}^{\infty}
  \pi_{j+\frac{m}2}^{SL(2,{\mathbb{R}})} \boxtimes \mathcal H^{j}(\mathbb{R}^{m}),
\end{align*}
where
 ${\mathcal{H}}^j({\mathbb{R}}^m)$ denotes
 the natural representation of $O(m)$
 (or $SO(m)$)
 on the space
 of harmonic polynomials on ${\mathbb{R}}^m$
 of degree $j$
 and $\pi_{b}^{SL(2,{\mathbb{R}})}$
 stands for the irreducible unitary lowest weight representation
 of $SL(2,{\mathbb{R}})$
 (or its covering group)
 with minimal $K$-type $b$.

These considerations bring us
 to interpolate operators
 occurring two minimal representations
 of $SO_0(m+1,2)\,\widetilde{}$ and $Sp(m,{\mathbb{R}})$.  
For this,
 we take $a >0$
 to be a deformation parameter,
 and define
\[
  \widetilde{h}_a:=\frac 2 a E + \frac {m+a-2}{a}, 
  \quad
  \widetilde{e}_a:= \frac {\sqrt{-1}}{a}|x|^a, 
  \quad
  \widetilde{f}_a:=\frac {\sqrt{-1}}{a}|x|^{2-a}\Delta.  
\]
The operators \eqref{eqn:Spsl} in the Weil representation
 corresponds to the case $a=2$, 
 and the operators \eqref{eqn:Osl} for $SO_0(m+1,2)\,\widetilde{}$
 corresponds to the case $a=1$.  
They extend to self-adjoint operators
 on the Hilbert space $L^2({\mathbb{R}}^m, |a|^{a-2}dx)$, 
 form an ${\mathfrak{sl}}_2$-triple, 
and lift to a unitary representation 
 of the universal covering group $SL(2,{\mathbb{R}})\,\widetilde{}
$
 of $SL(2,{\mathbb{R}})$
 for every $a>0$.  
The Hilbert space decomposes
 into a multiplicity-free discrete sum
 of irreducible unitary representations
 of $SL(2,{\mathbb{R}})\,\widetilde{} \times O(m)$
 as follows:
\[
  L^2({\mathbb{R}}^m, |x|^{a-2}dx)
  \simeq
  \sideset{}{^\oplus}{\sum}_{j=0}^{\infty}
  \pi_{\frac{2j+m-2}{a}+1}^{SL(2,{\mathbb{R}})} 
  \boxtimes {\mathcal{H}}^j({\mathbb{R}}^m).  
\]

The discrete decomposition
of ${\mathfrak {sl}}_2$-modules
 becomes a tool to generalize 
 the study of the unitary inversion operator
 ${\mathcal{F}}_{\Xi}$
 and the holomorphic semigroup 
 in \cite{xkmano1, xkmano2}
 to the following settings:
\begin{itemize}
\item[$\bullet$]
Dunkl operators
 (with Ben Sa\"id and {\O}rsted \cite{xbko}),
\item[$\bullet$]
Conformal group of
 Euclidean Jordan algebras
 (with Hilgert and M{\"o}llers \cite{xhkmo}).
\end{itemize}

\section{Quantization of Kostant--Sekiguchi Correspondence}

In this section we discuss Theorem \ref{thm:khmL2}
 in a special case
 where $V$ is Euclidean,
 equivalently,
 $G/K$ is a tube domain, 
 and explain 
 a recent work 
 \cite{xhkmo}
 with Hilgert,  
 M\"ollers, 
 and {\O}rsted 
 on the construction of a new model
 (a Fock-type model)
 of minimal representations
 with highest weights
 and a generalization of the classical Segal--Bargmann transform, 
 which we called
 a \lq{geometric quantization}\rq\
 of the Kostant--Sekiguchi correspondence.
In the underlying idea,
 the discretely decomposable
 restriction of ${\mathfrak {sl}}(2,{\mathbb{R}})$, 
 which appeared in \cite{xkmano1}, 
 plays again an important role.

We recall
 (e.g., \cite{xfolland, Howe})
 that the classical Fock space ${\mathcal{F}}({\mathbb{C}}^m)$
 is a Hilbert space
 in the space ${\mathcal{O}}({\mathbb{C}}^m)$
 of holomorphic functions defined by
\[
  {\mathcal{F}}({\mathbb{C}}^m):=
  \{f \in {\mathcal{O}}({\mathbb{C}}^m)
    :\int_{{\mathbb{C}}^m} |f(z)|^2 e^{-|z|^2} dz < \infty
   \},
\]
 and that the Segal--Bargmann transform
is a unitary operator
\[
  {\mathcal{B}}:
  L^2({\mathbb{R}}^m) \overset {\sim} \to {\mathcal{F}}({\mathbb{C}}^m),
  \quad
  u \mapsto ({\mathcal{B}}u)(z)
           :=\int_{{\mathbb{R}}^m} K_{\mathcal{B}}(x,z)f(x) d x,
\]
with the kernel
\[
  K_{\mathcal{B}}(x,z):=\exp(-\frac 1 2\langle z, z\rangle
                +2 \langle z, x\rangle
                -\langle x, x\rangle).
\]
{}From a representation theoretic viewpoint,
 the classical Segal--Bargmann transform
 intertwines the two models
 of the Weil representation
 of the metaplectic group $Mp(m,{\mathbb{R}})$,
 namely,
 the Schr{\"o}dinger model on $L^2({\mathbb{R}}^m)$
 and the Fock model on ${\mathcal{F}}({\mathbb{C}}^m)$.

In order to find a natural generalization 
 of this classical theory,
 we begin by examining 
 how one may rediscover the classical Fock model.
Our idea is to use 
 the action of ${\mathfrak {sl}}_2$, 
more precisely, 
 a \lq{holomorphically extended representation}\rq\
 of an open semigroup of $SL(2,{\mathbb{C}})$
 rather than a unitary representation 
 of $SL(2,{\mathbb{R}})$ itself.
For this, 
 we take a 
standard basis of ${\mathfrak {sl}}(2, {\mathbb{R}})$
 as 
\begin{equation}
\label{eqn:hef}
  h:= \begin{pmatrix} 1 & 0 \\ 0 & -1 \end{pmatrix},
  \quad
  e:= \begin{pmatrix} 0 & 1 \\ 0 & 0 \end{pmatrix},
  \quad
  f:= \begin{pmatrix} 0 & 0 \\ 1 & 0 \end{pmatrix}.  
\end{equation}
They satisfy the following Lie bracket relations:
$
  [h,e]=2e, \quad [h,f]=-2f, \quad [e,f]=h.
$
We set
\begin{align}
    k:=& i(-e+f)
      =\begin{pmatrix} 0 & -i \\ i & 0 \end{pmatrix}, 
\notag
\\
\label{eqn:BFc}
c_1:=&
     \begin{pmatrix}
     1 & -i
     \\
     -\frac{i}{2} & \frac{1}{2}
     \end{pmatrix}
=
     \begin{pmatrix}
     2i & 0
     \\
     0 &\frac{1}{2i}
     \end{pmatrix}
     \begin{pmatrix}
     1 & -\frac{i}{2}
     \\
     0 & 1
     \end{pmatrix}
     \begin{pmatrix}
     0 & -1
     \\
     1 &0
     \end{pmatrix}
     \begin{pmatrix}
     1 & i
     \\
     0 & 1
     \end{pmatrix}.
\end{align}

By a simple matrix computation 
 we have:
\begin{equation}
\label{eqn:w}
  \exp(-\frac{t}{2}k)|_{t=i \pi}
  = \begin{pmatrix} 0 & -1 \\ 1 & 0 \end{pmatrix}
  \in SL(2,{\mathbb{R}}).
\end{equation}
The formula $\operatorname{Ad}(c_1) k = h$ shows
 that $c_1 \in SL(2,{\mathbb{C}})$ gives a Cayley transform.
Correspondingly,
 the Bargmann transform may be interpreted
as
\[
{\mathcal {B}}= \text{\lq} \pi\circ \iota (c_1) \text{\rq}.  
\]
The right-hand side is not well-defined.  
We need an analytic continuation
 in the Schr{\"o}dinger model
 and a lift in the diagram below:
\begin{align*}{}
  &SL(2,{\mathbb{R}})\,\widetilde{} \overset \iota \to G
   \overset {\pi}\to GL(L^2(\Xi))
\\
  &\hphantom{mii}\downarrow
\\
c_1 \in SL(2,{\mathbb{C}}) \supset &SL(2,{\mathbb{R}})
\end{align*}

To be more precise,
 we write $w \in G$
 for the lift of \eqref{eqn:w}
 via
 $d \iota: {\mathfrak {sl}}(2,{\mathbb{R}}) \hookrightarrow {\mathfrak {g}}$.
Since the action of the maximal parabolic subgroup $P$
 on $L^2(\Xi)$
 is given by the translation
 and the multiplication
 of functions,
 it is easy to see what \lq{$\pi(p)$}\rq\
 should look like for $p \in P_{\mathbb{C}}$.
Therefore,
 we could give an explicit formula 
 for the (generalized) Bargmann transform
$
  {\mathcal{B}}=\text{\lq{$\pi \circ \iota(c_1)$}\rq}
$
 if we know the closed formula
 of the unitary inversion:
\[
  {\mathcal{F}}_{\Xi}
 ={\mathcal{F}}(w)
 \equiv
 {\mathcal{F}} \circ \iota \begin{pmatrix} 0 & -1 \\ 1 & 0\end{pmatrix}.
\]
Of course,
 this is not a rigorous argument,
 and $\pi(p)$ does not leave $L^2(\Xi)$
 invariant.  
However, 
 the formula \eqref{eqn:BFc}
 suggests what the function space $\pi \circ \iota(c_1)(L^2(\Xi))$
 ought to be, 
 and led us to an appropriate generalization of 
the classical Fock space
 as follows:
\begin{equation}
\label{eqn:FO}
{\mathcal{F}}({\mathbb{O}}_{\operatorname{min}}^{K_{\mathbb{C}}})
:=
\{F \in {\mathcal{O}}({\mathbb{O}}_{\operatorname{min}}^{K_{\mathbb{C}}}):
  \int_{{\mathbb{O}}_{\operatorname{min}}^{K_{\mathbb{C}}}}
  |F(z)|^2 \widetilde K_{\lambda-1}(|z|) d\nu(z) < \infty
\}.
\end{equation}
Here ${\mathbb{O}}_{\operatorname{min}}^{K_{\mathbb{C}}}$
 is the minimal nilpotent $K_{\mathbb{C}}$-orbit
 in ${\mathfrak{p}}_{\mathbb{C}}$
 which is the counterpart
 of the minimal (real) nilpotent coadjoint orbit
 ${\mathbb{O}}_{\operatorname{min}}^{G_{\mathbb{R}}}$
 in ${\mathfrak {g}}^{\ast}\simeq {\mathfrak {g}}$
under the Kostant--Sekiguchi correspondence
 \cite{xSeJ},
 see Figure 8.1.
Thus the generalized Fock space
 ${\mathcal{F}}({\mathbb{O}}_{\operatorname{min}}^{K_{\mathbb{C}}})$
 is a Hilbert space
 consisting of $L^2$-holomorphic functions
 on the complex manifold 
 ${\mathbb{O}}_{\operatorname{min}}^{K_{\mathbb{C}}}$
 against the measure 
 given by a renormalized $K$-Bessel function 
 $\widetilde K_{\lambda-1}(|z|)d\nu(z)$
 (see the comments after \eqref{eqn:Ftube}).  
\begin{alignat*}{2}
&{\mathfrak {g}}
&&{\mathfrak {p}}_{\mathbb{C}}
\\
&\cup &&\cup
\\
&{\mathbb{O}}_{\operatorname{min}}^{G_{\mathbb{R}}}
\overset {\text{Kostant--Sekiguchi}} 
{\leftarrow
\hskip -0.35pc 
-
\hskip -0.35pc
-
\hskip -0.35pc
-
\hskip -0.35pc
-
\hskip -0.35pc
-
\hskip -0.35pc
\rightarrow}
&&{\mathbb{O}}_{\operatorname{min}}^{K_{\mathbb{C}}}
\\
&\cup{\small{\text{Lagrangian}}}
&&
\\
&\Xi
&&
\end{alignat*}
\vspace {-0.5cm}
\[\small{\text{Figure 8.1 (minimal nilpotent orbits in ${\mathfrak {g}}$
 and ${\mathfrak {p}}_{\mathbb{C}}$)}}
\]
\vskip 1pc

We recall that $\Xi$ is a Lagrangian submanifold
 of ${\mathbb{O}}_{\operatorname{min}}^{G_{\mathbb{R}}}$,
 and $K_{\mathbb{C}}$ acts holomorphically
 on ${\mathbb{O}}_{\operatorname{min}}^{K_{\mathbb{C}}}$.
Then as a \lq{quantization}\rq\
 of the Kostant--Sekiguchi correspondence, 
 we define 
 the generalized Bargmann transform 
$
  {\mathcal {B}}:L^2(\Xi) \to \mathcal{F}({\mathbb{O}}_{\operatorname{min}}^{K_{\mathbb{C}}})
$
 by 
\[
  f \mapsto \Gamma(\lambda)e^{-\frac 1 2 \operatorname{tr}(z)}
            \int_{\Xi} \widetilde{I}_{\lambda-1}(2 \sqrt{(z|x)})
            e^{-\operatorname{tr}(x)}f(x) d \mu(x),
\]
whereas the unitary inversion operator ${\mathcal {F}}_{\Xi}$
 is given by
\begin{equation}
\label{eqn:Ftube}
  ({\mathcal {F}}_{\Xi}f)(y)
 =2^{-r \lambda}\Gamma(\lambda)
        \int_{\Xi}\widetilde J_{\lambda -1}(2 \sqrt{(x|y)})
        f(x) d\mu(x).
\end{equation}
Here $r=\operatorname{rank} G/K$,
 $(\,\,|\,\,)$ denotes the trace form
 of the Jordan algebra $V$,
 and $\lambda = \frac 1 2 \dim _{\mathbb{R}}{\mathbb{F}}$
 if $V = \operatorname{Herm}(k,{\mathbb{F}})$
 with ${\mathbb{F}}={\mathbb{R}}, {\mathbb{C}}$,
 quaternion ${\mathbb{H}}$,
 or the octonion ${\mathbb{O}}$
 (and $k=3$)
 or $\lambda = \frac 1 2 (k-2)$
 if $V={\mathbb{R}}^{1,k-1}$.
$\widetilde J(t)$, $\widetilde I(t)$,
and $\widetilde K(t)$ are the renormalization of the
 $J$-, $I$-, and $K$-Bessel function,
 respectively, 
 following the convention of \cite{xkmanoAMS}.


\begin{acknowledgement}
The author warmly thanks Professor Vladimir Dobrev for his
hospitality during the ninth  International Workshop:
 Lie Theory and
 its Applications in Physics in Varna, Bulgaria, 20-26 June 2011.
Thanks are also due to anonymous referees
 for careful readings.  
\\
The author was partially supported by
        Grant-in-Aid for Scientific Research (B) (22340026), Japan
        Society for the Promotion of Sciences.
\end{acknowledgement}

\end{document}